\documentclass[12pt,a4paper]{amsart}
\usepackage{amsmath,amssymb,amsthm}
\usepackage{enumitem}
\usepackage{ifthen,verbatim}
\usepackage{mathrsfs}
\usepackage{color}
\usepackage{url}
\usepackage[english]{babel}
\usepackage{here}
\usepackage[T1]{fontenc}
\usepackage{floatflt,graphicx,graphics}
\usepackage{a4wide}
\usepackage{cite}
\usepackage{ifthen}
\newtheorem{theorem1}{Theorem}[section]
\newtheorem{lemma1}{Lemma}[section]
\newtheorem{remark1}{Remark}[section]
\nonstopmode
\begin{document}

\title[Asymptotics of the conformal modulus of a nonsymmetric]{Asymptotics of the conformal modulus of a nonsymmetric unbounded  doubly-connected domain under stretching}

\author{Giang~V.~Nguyen}
\address{Kazan Federal University, Kremlyovskaya str. 35, Tatarstan, 420008,
Russian Federation}
\email{nvgiang.math@gmail.com} 

\author{Semen R.~Nasyrov}
\address{Kazan Federal University, Kremlyovskaya str. 35, Tatarstan, 420008,
Russian Federation}
\email{snasyrov@kpfu.ru} 


\begin{abstract}

We describe the asymptotic behavior of the conformal modulus of an unbounded doubly-connected domain, non-symmetric with respect to the coordinate axes, when stretched in the direction of the abscissa axis with coefficient $H\to +\infty$. Therefore, we give a partial
answer to a problem, suggested by M.~Vourinen, for the case of
an arbitrary unbounded domain.

\end{abstract}



\maketitle

\section{Introduction}
\subsection{The conformal modulus of a doubly-connected domain}
Let $D$ be a doubly-con\-nec\-ted planar domain with nondegenerate boundary
components, then $D$ is conformally equivalent to an annulus $\{r< |z|< R\}$ (see, e.g., \cite[Chap.~V, \S~1]{goluzin1969geometric}, \cite[Sect.~3]{kuhnau2005conformal}), and its conformal modulus $m(D)$ is the number defined as follows:
$$
m(D):=\dfrac{1}{2\pi}\log\dfrac{R}{r}\,.
$$
Another approach to define the conformal modulus  is in using  the concept of extremal length of curve family (see, e.g., \cite[Chap.~1, Sect.~D]{ahlfors2006lectures}, \cite[Chap.~3,
Sect.~11]{kuhnau2005conformal}). Specifically,
$$
m(D)=\lambda(\Gamma)
$$
where $\lambda(\Gamma)$ is the extremal length of the family $\Gamma$ of curves in the domain $D$ which join the boundary components of
$D$. Moreover,
$$
m(D)={1}/{\lambda(\Gamma')}
$$
where $\Gamma'$ is the family of curves  in $D$ separating its boundary components.

Conformal moduli of doubly-connected domains are also closely related to the conformal moduli of quadrilaterals (see, e.g., \cite{ahlfors2006lectures,kuhnau2005conformal,papamichael2010numerical}). A quadrilateral
$(Q;z_1,z_2,z_3,z_4)$ is a Jordan domain $Q$ with four fixed
points (vertices) $z_k$, $1\le k\le 4$,  on its boundary; the increasing of the index $k$ corresponds to the positive bypass of the boundary $\partial Q$.

For a given quadrilateral
$\mbox{\boldmath$Q$}=(Q;z_1,z_2,z_3,z_4)$, there exists a
conformal mapping of $Q$ onto a quadrilateral $(\Pi;0, 1,1+hi,hi)$, where $\Pi=[0,1]\times[0,h]$,
$h>0$, keeping the vertices. By definition,  the conformal modulus
$m(\mbox{\boldmath$Q$})$ of $\mbox{\boldmath$Q$}$ is equal to $h$.
The conformal modulus $m(\mbox{\boldmath$Q$})$ coincides with the
extremal length of the family of curves connecting in $Q$ the boundary arcs
$z_1z_2$ and $z_3z_4$; it is also the inverse to the extremal
length of the family of curves in $Q$ joining the arcs $z_2z_3$
and $z_4z_1$. The quadrilateral
$\mbox{\boldmath$Q$}^*=(Q;z_2,z_3,z_4,z_1)$ is called conjugate to
$\mbox{\boldmath$Q$}$. It is clear that
$m(\mbox{\boldmath$Q$}^*)=m^{-1}(\mbox{\boldmath$Q$})$.
If $\mbox{\boldmath$Q$}^c=(Q^c;z_4,z_3,z_2,z_1)$, where $Q^c=\overline{\mathbb{C}}\setminus Q$, then the conformal modulus of $\mbox{\boldmath$Q$}^c$ is called the exterior modulus (see, e.g., \cite{vuorinen2013exterior}).

Over the years, there have been obtained many results on the
conformal moduli of doubly-connected domains, as well as on moduli of quadrilaterals  (see, e.g.,  \cite{dautova2018,dautova2019,dautova2021conformal, hakula2021conformal,hrv1,hrv2,kuhnau2005conformal,nass,nasser2021circular,nasyrov2015riemann,nvgiang2021,nasyrov2021moduli,vuorinen2013exterior}).

Conformal moduli of doubly-connected domains and quadrilaterals are conformally invariant. They are also quasiinvariant under quasiconformal mappings (see, e.g., \cite[Chap.~2, Sect.~A]{ahlfors2006lectures}, \cite[Chap.~3, Sect.~22]{kuhnau2005conformal}). If $f$ is an $H$-quasiconformal mapping, then
$$
\dfrac{1}{H}m(\Omega)\leq m(f(\Omega))\leq H m(\Omega)
$$
where $\Omega$ is a doubly-connected domain. The same is valid for conformal moduli of quadrilaterals.

Prof. M.~Vuorinen posed the problem of investigating the distortion of the conformal moduli of doubly-connected  domains under the action of the $H$-quasiconformal mappings,
\begin{equation}\label{cttuabg}
f_H:x+iy\mapsto Hx+iy\quad (H>1);
\end{equation}
that is the stretching along the abscissa axis with coefficient $H$. In particular, it is interesting to study its asymptotic behavior as $H\to+\infty$ (see, e.g., \cite{nasyrov2015riemann}).

In \cite{dautova2018}, \cite{dautova2019}, \cite{nasyrov2015riemann}, \cite{nvgiang2021} the results on solving the Vuorinen problem were obtained for the cases where $D$ is a rectangular frame, a symmetric or non-symmetric bounded doubly-connected domain  or unbounded doubly-connected domains, symmetric with respect to the coordinate axes. Here we study the problem by Vuorinen for the case of an unbounded doubly-connected domain of sufficiently arbitrary form, non-symmetric with respect to the coordinate axes.
\medskip

\subsection{The main result}

Let $\mathfrak{S}$ be the set of all domains $\Omega$ which satisfy the following conditions:

(i) $\Omega$ is an unbounded doubly-connected domain;

\begin{figure}[ht] \centering
\includegraphics[width=3.5 in]{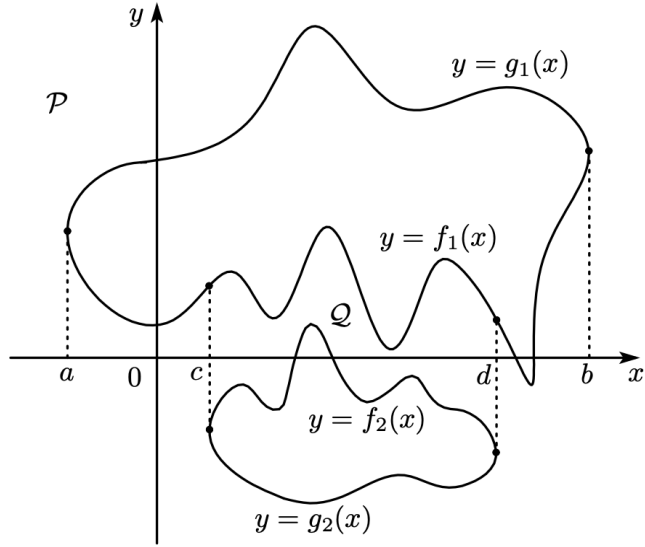}
\caption{A domain $\Omega$ from the set $\mathfrak{S}$.}\label{f6}
\end{figure}

(ii) the boundary of $\Omega$ consists of the graphs of continuous functions:  functions  $g_1$, $f_1$ are given on a segment $[a, b]$ and $f_1(x)\le g_1(x)$, $x\in[a, b]$; the functions  $g_2$, $f_2$ are given on a segment $[c, d]\subset[a,b]$ and $g_2(x)\le f_2(x)<f_1(x)$, $x\in[c, d]$; moreover, $f_1(a)=g_1(a), f_1(b)=g_1(b)$; $f_2(c)=g_2(c), f_2(d)=g_2(d)$ (see. Fig.~\ref{f6}).

We will describe the behavior of the conformal modulus of the domain $\Omega_H$ obtained from  $\Omega\in\mathfrak{S}$ by the stretching $f_H$, given by \eqref{cttuabg}, as $H\to +\infty$.

\begin{theorem1}\label{main}
If $\Omega\in\mathfrak{S}$ and $\Omega_H=f_H(\Omega)$, then
\begin{equation}\label{ctulmd1}
m(\Omega_H)\sim \dfrac{1}{\gamma H}\,, \quad H\to\infty, \quad
\mbox{\rm where }\ \gamma
=\int\limits_c^d\dfrac{dx}{f_1(x)-f_2(x)}\,.
\end{equation}
Moreover,
\begin{equation}\label{bdtmodulOmegaH}
m(\Omega_H)\leq\dfrac{1}{\gamma H}\quad  \mbox{\rm for every}\quad
H>0.
\end{equation}
\end{theorem1}

\section{Local versions of the Rad\'o convergence theorem}

In this section we will obtain a local version of the Rad\'o convergence theorem and its generalization for a sequence of doubly-connected domains with degenerate moduli (see, e.g., \cite[Theorem~3]{dautova2019}); these results will be used in the proof of Theorem~\ref{main}.

As it is well known, every conformal mapping $f$ between Jordan domains can be extended, by continuity, up to a homeomorphism of the closures of the domains (see, e.g., \cite[Chap.~2, Sect.~3]{goluzin1969geometric}); for simplicity of notation, we will also denote this homeomorphism by~$f$.

Firstly, we recall the classical Rad\'o's theorem on uniform convergence of conformal mappings. Let $E$ denote the unit disk $\{|z|<1\}$ (see, e.g., \cite[Chap.~2, Sect.~5]{goluzin1969geometric}).

\begin{theorem1}[Rad\'o]\label{thr_rado_ge}
Let $D_n$ be a sequence of Jordan domains, bounded
by curves $\Lambda_n$, which converges to a kernel with respect to a point $z_0$, and let this kernel be a Jordan domain $D$ bounded by a curve $\Lambda$. Then a sequence of continuous functions $f_n : \overline{E} \to \overline{D}_n$, mapping conformally ${E}$ onto $D_n$ and normalized by the conditions $f_n(0) = z_0$, $f_n' (0) > 0$, converges uniformly to the continuous function $f\,:\,\overline{E} \to \overline{D}$, $f(0) = z_0, f'(0) > 0$ mapping conformally $E$ onto $D$ if and only if there exist homeomorphisms between the curves $\Lambda_n$ and $\Lambda$ such that, for any $\varepsilon > 0$, there is an $N$ such that the distance between the respective points of the curves $\Lambda_n$ and $\Lambda$ is smaller than $\varepsilon$ whenever $n\geq N$.
\end{theorem1}

\begin{remark1}\label{rem1}
It is easy to prove that, under the assumptions of Theorem~\ref{thr_rado_ge}, 
$f^{-1}_n\rightrightarrows f^{-1}$ on every compact $K$ such that $K\subset \overline{D}_n$ for sufficiently large $n$. In particular, $f^{-1}_n\rightrightarrows f^{-1}$ on $\bigcap_{n=1}^\infty \overline{D}_n$.
\end{remark1}

First we will establish a local version of Theorem~\ref{thr_rado_ge}.

\begin{figure}[ht] \centering
\includegraphics[width=4.7 in]{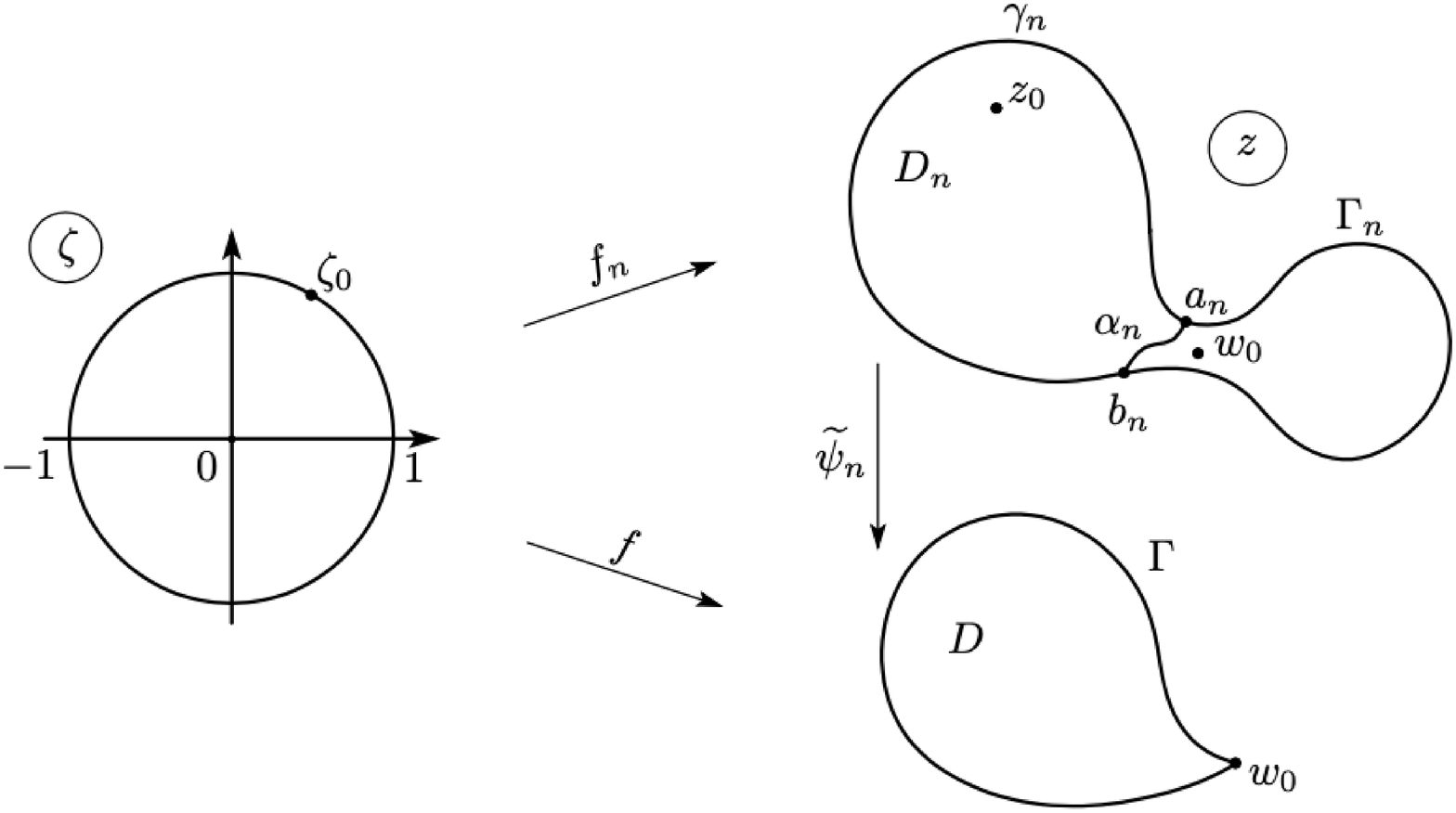}
\caption{The sequence of domains $D_n$ and its kernel.}\label{f1}
\end{figure}


\begin{theorem1}\label{thr_local_rado}
Assume that $D_n$ is a sequence of Jordan domains, bounded by curves $\Gamma_n$, which converges to a kernel with respect to a point $z_0$, and this kernel is a Jordan domain $D$ bounded by a curve $\Gamma$. Let $f_n:\overline{E}\to \overline{D}_n$ be a sequence of continuous functions, mapping conformally ${E}$ onto $D_n$ and normalized by the conditions $f_n(0) = z_0$, $
f_n'(0)>0,$ and $f:\overline{E}\to \overline{D}$ be a continuous function mapping conformally  ${E}$ onto $D$ and  satisfying the conditions $f(0) = z_0$, $f'(0)>0$. Let there exist a sequence of subarcs  $\gamma_n\subset \Gamma_n$ with endpoints $a_n$, $b_n$,  a sequence of continuous mappings $\psi_n: \gamma_n\to \Gamma$ and a point $w_0\in \Gamma$ such that:

$(i)$\ \textit{for every $n\ge 1$, $\psi_n$ is a homeomorphism of $\gamma_n\setminus\{a_n,b_n\}$ onto} $\Gamma\setminus \{w_0\}$;

$(ii)$\ for every $n\ge 1$, $\psi_n(a_n)=\psi_n(b_n)=w_0$;

$(iii)$\ for every $\varepsilon> 0$ there exists a number $N\geq 0$ such that
\begin{equation*}
\left|\psi_n(z)-z\right|<\varepsilon\quad \forall n\geq N\  \forall z\in\gamma_n;
\end{equation*}

$(iv)$\  there exists a sequence of crosscuts $\alpha_n$ of domains $D_n$ with endpoints $a_n$ and $b_n$, diameters of which tend to $0$ as $n\to\infty$.

Then $f_n$ converges locally uniformly to $f$ on $\overline{E}\setminus\{\zeta_0\}$, $\zeta_0=f^{-1}(w_0)$, and $f^{-1}_n$ converges uniformly to $f^{-1}$ on every compact subset of $\bigcap_{n=1}^\infty \overline{D}_n\setminus \{w_0\}$.
\end{theorem1}
\begin{proof}
I) It is clear that the sequence of Jordan domains $D'_n$, bounded with $\gamma_n\cup\alpha_n$, converges to $D$ as to a kernel with respect to $z_0$ (Fig.~\ref{f1}). Without loss of generality, we can assume that  $z_0\in D_n$, $n\ge 1$. Consider a sequence of conformal mappings
$$
\widetilde {f}_n: E \to D'_n
$$
with normalization $\widetilde {f}_n (0) = z_0 $, $\widetilde {f}_n '( 0)> 0 $.

\begin{lemma1}\label{wfn}
The sequence $\widetilde{f}_n$ converges uniformly to $f$ on the closed unit disk~$\overline{E}$.
\end{lemma1}

\begin{proof}
Using the conditions of the theorem, we can easily show that there exist homeomorphisms
$$
\widetilde {\psi}_n: \alpha_n \cup \gamma_n \to \Gamma,
$$
satisfying the Rad\'o condition, i.e. such that for any $\varepsilon> 0$ there exists $N\geq 0$ such that
\begin{equation*}
|\widetilde{\psi}_n(z)-z|<\varepsilon, \quad z\in\alpha_n \cup\gamma_n,
\end{equation*}
whenever $n\geq N$. By Theorem~\ref{thr_rado_ge}, we conclude that $\widetilde{f}_n$ converges uniformly to $f$ on~$\overline{E}$.
\end{proof}
\medskip

II) Here we will show that the preimages of the arcs $\alpha_n$ under the mappings $\widetilde {f}_n$ are arcs on the unit circle which are contracted to a point as $n\to\infty$.

According to the principle of the boundary correspondence under conformal mappings, every $\beta_n:=\widetilde{f}_n^{-1}(\alpha_n)$ is an arc on the unit circle, with the end points $\zeta_{n1}^a:=\widetilde{f}_n^{-1}(a_n)$ and $\zeta_{n1}^b:=\widetilde{f}_n^{-1}(b_n)$. Similarly, the points $\zeta^a_{n2}:=f^{-1}_n(a_n)$ and $\zeta^b_{n2}:=f^{-1}_n(b_n)$ are also on the boundaries of the domains $f^{-1}_n(D'_n)$.

\begin{lemma1}\label{lem1}
The diameters of $\beta_n$ tend to $0$ as $n\to\infty$. Moreover, these arcs are contracted to the  point $\zeta_0$, as $n\to\infty$.
\end{lemma1}

\begin{proof}
First we recall that by (iv) the diameters of the arcs $\alpha_n$ tend to 0, and by (ii) and (iii) their endpoints $a_n$ and $b_n$ converge to $w_0$. Therefore, $\alpha_n$ are contracted to the point $w_0$, as $n\to\infty$. Assume that $\mbox{\rm diam}(\beta_n)\not\to 0$, as $n\to\infty$.  Then there exists a subsequence $\beta_{n_k}$  of the sequence $\beta_{n}$ such that the  $\beta_{n_k}$ tends to some nondegenerate subarc $\beta$ of the unit circle. By Lemma~\ref{wfn}, we have  $\widetilde{f}_{n_k}\rightrightarrows f$ on $\overline{E}$, therefore,  $f\equiv w_0$ on $\beta$. By the boundary uniqueness theorem (see, e.g., \cite[Chap.~2, Sect.~3, Lemma~3]{goluzin1969geometric}) this implies that  $f\equiv w_0$  on the unit disk. The contradiction proves that $\mbox{\rm diam}(\beta_n)\to 0$, as $n\to\infty$.

If $\beta_n$ are not contracted to the point $\zeta_0$, then there exist a subsequence $\beta_{n_k}$ of $\beta_n$ converging to some $\zeta_1\neq \zeta_0$. But $\widetilde{f}_{n_k}$ converges uniformly to $f$ on the closed unit disk. For every sequence $\zeta_{n_k}\in\beta_{n_k}$ we have $\zeta_{n_k}\to \zeta_1$ and $f_{n_k}(\zeta_{n_k})\to w_0$, $k\to \infty$, since $f_{n_k}(\zeta_{n_k})\in \alpha_{n_k}$. But, because of uniform convergence,  $f(\zeta_1)=\lim_{k\to\infty} f_{n_k}(\zeta_{n_k})=w_0$, therefore, $\zeta_1=\zeta_0$. The contradiction proves the lemma.
\end{proof}

III) Consider the arcs $\sigma_n:=f_n^{-1}(\alpha_n)$. Denote by $\widetilde{D}_n$ the extended complex plane with the excluded arc $\beta_n$. Let $\widetilde{\sigma}_n$ be the arc symmetric to $\sigma_n$ with respect to the unit circle. Let $\widetilde{G}_n$ be the unbounded domain in the extended complex plane with the boundary $\sigma_n\cup\widetilde{\sigma}_n$.

\begin{lemma1}\label{lem4}
The sequences of domains $\widetilde{D}_n$ and $\widetilde{G}_n$ converge to the domain $\overline{\mathbb{C}}\setminus\{\zeta_0\}$, as to a kernel, as $n\to\infty$. Moreover, for every $\delta>0$ the sequences  $f_n^{- 1}\circ\widetilde{f}_n$ and $\widetilde{f}_n^{-1}\circ{f}_n$ converge locally uniformly to the identity mapping on the set $\overline{E}\setminus \{\zeta_0\}$.
\end{lemma1}

\begin{proof}
Since, by Lemma~\ref{lem1}, the arcs $\beta_n$ are contracted into the point $\zeta_0$, as $n\to+\infty$, we see that $\widetilde{D}_n$ converges to $\overline{\mathbb{C}}\setminus\{\zeta_0\}$.

Denote $g_n:=f_n^{- 1}\circ \widetilde{f}_n$.
By the symmetry principle \cite[Chap.~X, Sect.~1]{lang2003complex},  we can extend the mapping $g_n$ through the arc of the unit disk, complement to  $\beta_n$, in $\widetilde{D}_n$. The extended mapping, which we also, for simplicity of notation, denote by $g_n$, maps $\widetilde{D}_n$ onto $\widetilde{G}_n$. On the other hand, the family of the functions $g_n$ is a normal family in $\overline{\mathbb{C}}\setminus\{\zeta_0\}$. Actually, $g_n$ is bounded in $E$, therefore, in $E$ it forms a normal family. This follows that in  $\overline{\mathbb{C}}\setminus\overline{E}$ it is also normal. Thus, we only need to prove that $g_n$ is normal in a small neighborhood of every point $\omega\in\partial E$ distinct from $\zeta_0$. Let us choose a neighborhood $U$ of $\omega$ not containing the points $\zeta_0$, $0$ and $\infty$. Then in $U$ every $g_n$ does not take  the values $0$, $\infty$ and some $e^{i\phi_n}\in \partial E$. Consequently, $e^{-i\phi_n}g_n$ does not take three fixed values $0$, $1$ and $\infty$. By Montel's theorem (see, e.g., \cite[Chap.~2, Sect.~7]{goluzin1969geometric}, \cite[Chap.~7, Sect.~2]{remmert2013classical}, \cite[Chap.~2, Sect.~2]{schiff1993normal}), $e^{-i\phi_n}g_n$ is a normal family in $U$ and there exists a subsequence $e^{-i\phi_{n_k}}g_{n_k}$ converging locally uniformly in $U$ to some limit function $h$. But $\partial E$ is a compact set, therefore, passing, if necessary, to a subsequence, we can assume that the sequence $e^{i\phi_n}$ converges to some $e^{i\phi_0}\in\partial E$ in $U$. Therefore, $g_{n_k}$ converges to $e^{i\phi_0}h$. This means that the family $g_n$ is normal.

Consider any subsequence $g_{n_k}$ converging locally uniformly to some $g$. Since $g_{n_k}(0)=0$, we see that $g(0)=0$. If $g\equiv \mbox{\rm const}$, then $g\equiv 0$, but from the symmetry of $g_n$ we see that if $g_{n_k}(z)\to 0$, then $g_{n_k}(1/\overline{z})\to \infty$. The controversy proves that $g\not\equiv \mbox{\rm const}$. Then $g$ is a univalent function in $\overline{\mathbb{C}}\setminus\{\zeta_0\}$ as a nonconstant limit of univalent functions (see, e.g., \cite[Chap.~1, Sect.~2]{goluzin1969geometric}), and  $g(0)=0$, $g'(0)>0$. Moreover, $g(1/\overline{z})=1/\overline{g(z)}$. This conditions imply that $g(z)\equiv z$, $z\in \overline{\mathbb{C}}\setminus\{\zeta_0\}$. Therefore, by the generalized Carath\'eodory convergence theorem (see, e.g., \cite[Chap.~5, Sect.~5]{goluzin1969geometric}), $\widetilde{G}_n$ converges to  $\overline{\mathbb{C}}\setminus\{\zeta_0\}$ as $n\to\infty$. By the same theorem, we deduce that the sequence of inverse functions $\widetilde{f}_n^{- 1}\circ {f}_n$ converges locally uniformly in $\overline{\mathbb{C}}\setminus\{\zeta_0\}$ to the identity mapping.
\end{proof}

Now we will complete the proof of Theorem~\ref{thr_local_rado}.

From Lemmas~\ref{wfn} and \ref{lem4} we obtain that $f_n=\widetilde{f}_n\circ (\widetilde{f}_n^{- 1}\circ {f}_n)$ converges locally uniformly to $f$ in $\overline{E}\setminus \{\zeta_0\}$.
Now we will prove that $f^{-1}_n$ converges uniformly to $f^{-1}$ on every compact subset $K$ of $\bigcap_{n=1}^\infty \overline{D}_n\setminus \{w_0\}$. Assume the opposite. Then there exist $\varepsilon>0$ and a subsequence of points $z_{n_k}\in K$  such that
\begin{equation}\label{f1nk}
|f^{-1}_{n_k}(z_{n_k})-f^{-1}(z_{n_k})|\ge \varepsilon, \quad  k\ge 1.
\end{equation}
Denote $\zeta_{n_k}=f^{-1}(z_{n_k})$. Since $\zeta_{n_k}\in \overline{E}$ and $\overline{E}$ is compact, we may assume that $\zeta_{n_k}$ converges to some $\zeta^*\in \overline{E}$. Then  $z_{n_k}=f(\zeta_{n_k})\to f(\zeta^*)=:z^*$. It is cleat that $z^*\in K$, therefore, $z^*\neq w_0$. Consequently,   $\zeta^*\neq \zeta_0$. From the uniform convergence of $f^{-1}_{n_k}$ to
$f^{-1}$ on $K$ it follows that
$$f^{-1}_{n_k}(z_{n_k})=f^{-1}_{n_k}\circ f(\zeta_{n_k})\to \zeta^*.$$
Thus, we have $f^{-1}_{n_k}(z_{n_k})\to f^{-1}(z^*)=\zeta^*$, as $k\to\infty$, but this contradicts to \eqref{f1nk}.

Therefore, Theorem \ref{thr_local_rado} is fully verified.
\end{proof}
\medskip

Now we will state a theorem which is a local version of Rad\'o type theorem proven in \cite[Theorem~3]{dautova2019}. In the cited theorem the case of convergence of univalent functions defined on a sequence of doubly-connected domains   converging to a simply-connected domain $G$ is considered. Here we investigate the case where the uniform convergence is violated at some point.

Let $\mathcal{D}_n$ be a sequence of non-concentric annuli, each of which is bounded by the unit circle $\{|\zeta|=1\}$ and the circle $\{|\zeta+ir_n/2|=r_n/2\}$, where $0<r_n<1$ and $\lim\limits_{n\to\infty}r_n=1$. Then, $\mathcal{D}_n$ converges to the circular lune $\mathcal{D}$ bounded by the unit circle and the circle $\{|\zeta+i/2|=1/2\}$.

\begin{figure}[ht] \centering
\includegraphics[width=4.5 in]{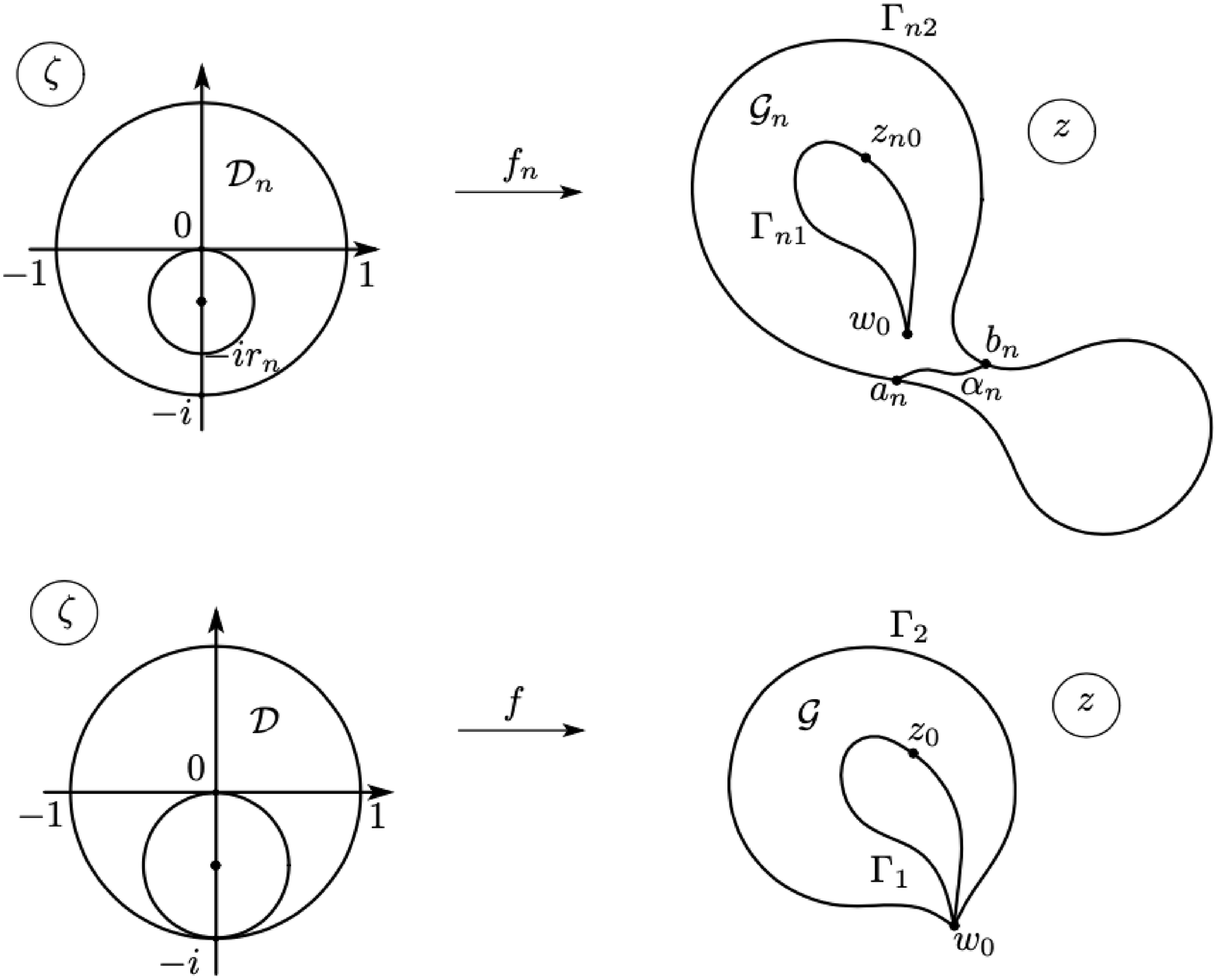}
\caption{The sequence of mappings $f_n$ and its limit.}\label{f3}
\end{figure}

Let $\mathcal{G}_n$ be a sequence of doubly-connected domains in $\mathbb{C}$, each of which is bounded by Jordan curves $\Gamma_{n1}$ and  $\Gamma_{n2}$, and $\Gamma_{n2}$ is the exterior component of $\partial \mathcal{G}_n$. Let there exist a sequence of crosscuts $\alpha_n$ of domains $\mathcal{G}_n$ with endpoints $a_n$, $b_n$, such that $\lim\limits_{n\to\infty}a_n= \lim\limits_{n\to\infty}b_n=w_0$ for some $w_0\in \mathbb{C}$ and the diameters of $\alpha_n$ tend to $0$, as $n\to\infty$. For every $n$ we denote by $\gamma_n$ a subarc of $\Gamma_{n2}$ with endpoints $a_n, b_n$ such that $\gamma_n\cup\alpha_n$ is the boundary of the Jordan domain containing $\Gamma_{n1}$.  Let the sequence $\mathcal{G}_n$ converge to a simply-connected domain $\mathcal{G}$, the boundary of which consists of two Jordan curves $\Gamma_1$ and $\Gamma_2$, and  $\Gamma_1\cap\Gamma_2=\{w_0\}$ (Fig.~\ref{f3}).

\begin{theorem1}\label{thr_rado_ge1}
Let $f_n$ be a sequence of functions mapping conformally $\mathcal{D}_n$ onto $\mathcal{G}_n$, and $f_n(0)=z_{n0}\in\Gamma_{n1}$. Let the function $f$ map $\mathcal{D}$ onto $\mathcal{G}$, and $f(0)=z_0\in\Gamma_1$, where $z_0\ne w_0$. Assume that $\lim\limits_{n\to\infty}z_{n0}=z_0$, and there exist homeomorphisms
$$\psi_{n1}:\Gamma_1\to\Gamma_{n1},\,\,\, \psi_{n2}:\Gamma_2\to\Gamma_{n2}$$
such that
\begin{equation}\label{eq_rado_ge1}
\forall\varepsilon>0\ \ \exists N:\quad\forall n\geq N\quad |\psi_{n1}(t)-t|<\varepsilon,\ t\in\Gamma_1;\qquad |\psi_{n2}(s)-s|<\varepsilon,\ s\in\Gamma_2.
\end{equation}
Then, $f_n$ converges locally uniformly to $f$ on $\overline{\mathcal{D}}\setminus\{-i\}$ and $f^{-1}_n$ converges uniformly to $f^{-1}$ on every compact subset of $\bigcap_{n=1}^\infty \overline{\mathcal{G}}_n\setminus \{w_0\}$.
\end{theorem1}

To prove the theorem, we need to establish some auxiliary results.

Denote $\alpha_n':=f_n^{-1}(\alpha_n), a_n':=f_n^{-1}(a_n), b_n':=f_n^{-1}(b_n)$. According to the principle of the boundary correspondence under conformal mappings, $\alpha_n'$ is a crosscut in $\mathcal{D}_n$, with endpoints $a_n', b_n'$ on the unit circle.

\begin{lemma1}\label{lem1a}
The crosscuts $\alpha'_n$ are contracted to the point $(-i)$, as $n\to\infty$.
\end{lemma1}

The proof of Lemma~\ref{lem1a} is similar to that of Lemma~\ref{lem1}.\medskip

\begin{figure}[ht] \centering
\includegraphics[width=4.5 in]{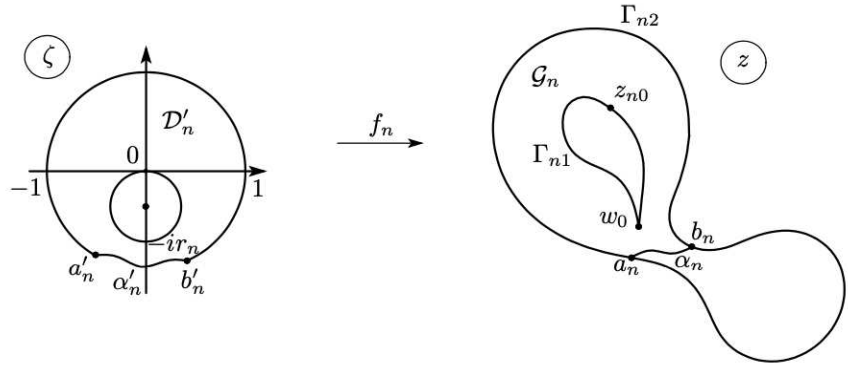}
\caption{The sequence of mappings $f_n$ on subdomains $\mathcal{D}'_n$.}\label{f4}
\end{figure}

Let $\mathcal{D}'_n$ be the sequence of doubly-connected domains, each of which is bounded by the circle $\{|\zeta+ir_n/2|=r_n/2\}$ and $\alpha'_n\cup\alpha''_n$, where $\alpha''_n$ is the arc of the unit circle $\{|\zeta|=1\}$ with endpoints $a'_n$, $b'_n$ such that the curve $\alpha'_n\cup\alpha''_n$ contains the circle $\{|\zeta+ir_n/2|=r_n/2\}$ in its interior (Fig.~\ref{f4}).

\begin{figure}[ht] \centering
\includegraphics[width=4.5 in]{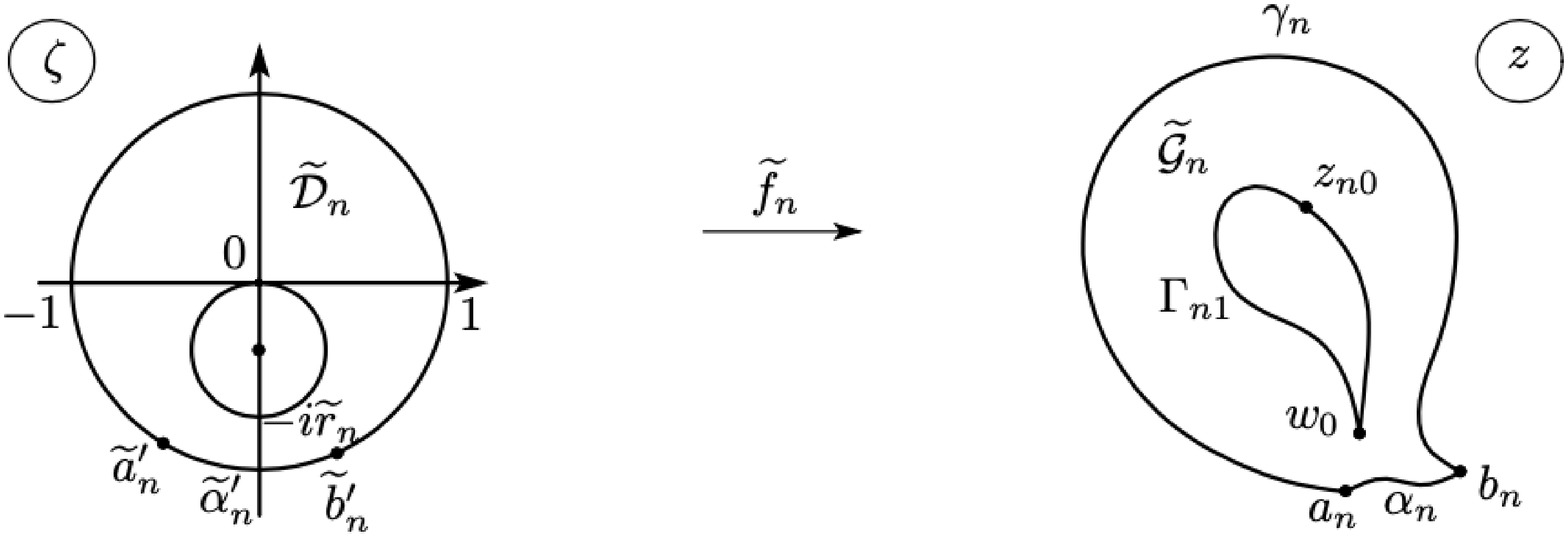}
\caption{The sequence of mappings $\widetilde{f}_n$.}\label{f5}
\end{figure}

Let $\widetilde{\mathcal{G}}_n$ be a part of $\mathcal{G}_n$, bounded by $\gamma_n\cup\alpha_n$ and $\Gamma_{n1}$ (Fig.~\ref{f5}). From Lemma \ref{lem1a} and the conditions on $\alpha_n $, we easily deduce the following lemma.

\begin{lemma1}
The sequences of domains $\mathcal{D}_n'$ and $\widetilde{\mathcal{G}}_n$ converge to the domains $\mathcal{D}$ and $\mathcal{G}$, respectively, as $n\to\infty$.
\end{lemma1}

Denote by $\widetilde{f}_n$ a sequence of functions mapping conformally  $\widetilde{\mathcal{D}}_n$ onto $\widetilde{\mathcal{G}}_n$ such that $\widetilde{f}_n(0)=z_{n0}\in\Gamma_{n1}$; here $\widetilde{\mathcal{D}}_n$ are nonconcentric annuli bounded by the unit circle and the circles $\{|\zeta+i\widetilde{r}_n/2|=\widetilde{r}_n/2\}$, where $0<\widetilde{r}_n<1$ and $\lim\limits_{n\to\infty}\widetilde{r}_n=1$. The arc $\widetilde{\alpha}_n':=\widetilde{f}^{-1}_n(\alpha_n)$ is a subarc of the unit circle with endpoints $\widetilde{a}_n':=\widetilde{f}_n^{-1}(a_n), \widetilde{b}_n':=\widetilde{f}_n^{-1}(b_n)$. Applying the same arguments as in the proof of Lemma \ref{lem1}, we obtain that the arcs $\widetilde{\alpha}'_n$ are contracted to the point $(-i)$, as $n\to\infty$.

It is easy to see that the sequence $\widetilde{\mathcal{D}}_n$ tends to $\mathcal{D}$ as $n\to\infty$.

Let the conditions of Theorem \ref{thr_rado_ge1} hold. Then it is evident that there exist homeomorphisms $\widetilde{\psi}_{n1}:\Gamma_1\to\Gamma_{n1}$ and $\widetilde{\psi}_{n2}:\Gamma_2\to\gamma_n\cup\alpha_n$ such that
\begin{equation}\label{eq_lem_rado}
\forall\varepsilon>0\ \exists N:\quad\forall n\geq N\quad |\widetilde{\psi}_{n1}(t)-t|<\varepsilon,\ t\in\Gamma_1;\quad |\widetilde{\psi}_{n2}(s)-s|<\varepsilon,\ s\in\Gamma_2.
\end{equation}
From this fact and the generalized Rad\'o theorem \cite[Sect.~4, Theorem.~6]{dautova2018} we obtain

\begin{lemma1}\label{lem_rado}
The sequence $\widetilde{f}_n$ converges uniformly to $f$ in $\overline{\mathcal{D}}$.
\end{lemma1}

Denote by ${\mathcal{D}}'_{n1}$ and $\widetilde{\mathcal{D}}_{n1}$ the domains which obtained from the domains $\mathcal{D}'_n$ and $\widetilde{\mathcal{D}}_n$ by gluing them with their reflections with respect to the arcs on the unit circle $|\zeta|=1$ with endpoints $a'_n, b'_n$ and $\widetilde{a}'_n, \widetilde{b}'_n$, complement to $\alpha_n$ and $\alpha'_n$ respectively. It is clear that the sequences
${\mathcal{D}}'_{n1}$ and $\widetilde{\mathcal{D}}_{n1}$  tend to the domain $\mathcal{D}_0=\{|z+1/2|>1/2, \,\text{Im} z>-1\}$ as $n\to \infty$.

Every function $f_n^{-1}\circ \widetilde{f}_n$ maps $\widetilde{\mathcal{D}}_n$ onto $\mathcal{D}'_n$.  According to the Riemann-Schwarz symmetry principle, the functions  $f_n^{-1}\circ \widetilde{f}_n$ can be extended to the domains $\widetilde{\mathcal{D}}_{n1}$; for convenience, we will also denote the extended functions by $f_n^{-1}\circ \widetilde{f}_n$; each of them  maps $\widetilde{\mathcal{D}}_{n1}$ onto
$\widetilde{\mathcal{D}}'_{n1}$.

Applying the same arguments as in the proof of Lemma \ref{lem4} we obtain

\begin{lemma1}\label{lem_rado1}
The sequence of functions $f_n^{-1}\circ\widetilde{f}_n$ converges locally uniformly to the identity mapping on $\mathcal{D}_0\setminus\{-i\}$.
\end{lemma1}

The same statement also holds for the inverse functions $\widetilde{f}_n^{-1}\circ f_n$. Then, with the use of Lemma~\ref{lem_rado}, we obtain that $f_n=\widetilde{f}_n\circ (\widetilde{f}_n^{- 1}\circ {f}_n)$ converges locally uniformly to $f$ in $\mathcal{D}_0\setminus\{-i\}$. The uniform convergence of $f_n^{-1}$ to $f^{-1}$ can be proved just in the proof of similar result from Theorem~\ref{thr_local_rado}.   This proves Theorem~\ref{thr_rado_ge1}.

\begin{remark1}\label{remfaminf}
Further we can apply Theorem~\ref{thr_rado_ge1} not only for sequences of conformal mappings but also for one parameter families of such mappings. We can also apply this theorem for the case of unbounded domains; then we understand the convergence in the spherical metric.
\end{remark1}
\section{Proofs of Theorem~\ref{main}}

First we will show that \eqref{bdtmodulOmegaH} holds.

Consider four the points $A(c, f_1(c)), B(c, f_2(c)), C(d, f_2(d)), D(d, f_1(d))$ on $\partial \Omega$. The vertical segments $AB$ and $CD$ separate $\Omega$ into two domains $\mathcal{Q}$ and $\mathcal{P}$. Let $\mathcal{Q}$ denote the bounded domain and $\mathcal{P}$ be unbounded, i.e.
$$\mathcal{Q}:=\{(x, y)\in\mathbb{C}\,:\, c\leq x\leq d, f_2(x)\leq y\leq f_1(x)\},$$
and $\mathcal{P}=\Omega\setminus \overline{\mathcal{Q}}$.

Consider the quadrilaterals
$$\mbox{\boldmath$Q$}:=\left(\mathcal{Q}; A, B, C, D\right),\quad \mbox{\boldmath$P$}:=\left(\mathcal{P}; D, C, B, A\right)$$
and their images $$
\mbox{\boldmath$Q$}_H:=f_H(\mbox{\boldmath$Q$})=\left(\mathcal{Q}_H; A_H, B_H, C_H, D_H\right),\quad  \mbox{\boldmath$P$}_H:=f_H(\mbox{\boldmath$P$})=\left(\mathcal{P}_H; D_H, C_H, B_H, A_H\right)
$$
under the mapping $f_H$. By the Gr\"otzsch inequality (see, e.g., \cite[Chap.~4, Sect.~6]{goluzin1969geometric},  \cite[Chap.~3, Sect.~7]{kuhnau2005conformal}),  we have
\begin{align}\label{m1}
m^{-1}(\Omega_H)\geq m(\mbox{\boldmath$Q$}_H) +m(\mbox{\boldmath$P$}_H).
\end{align}

From~\eqref{m1} we deduce that $m^{-1}(\Omega_H)\geq m(\mbox{\boldmath$Q$}_H)$ and from \cite[Theorem~4.2]{nvgiang2021}  it follows that
\begin{equation}\label{minv}
m(\mbox{\boldmath$Q$}_H)\ge \gamma H, \quad \gamma=\int\limits_{c
}^{d}\dfrac{dx}{f_1\left({x}\right)-f_2\left({x}\right)}.
\end{equation}
Therefore,
$m^{-1}(\Omega_H)\geq \gamma H$ and this implies~\eqref{bdtmodulOmegaH}. In particular, we have
\begin{equation}\label{momega}
m(\Omega_H)\to 0, \quad\text{ as }\quad H\to\infty.
\end{equation}

Denote by $\varphi_H$ a conformal mapping of $\Omega_H$ onto a concentric annulus
$$
\mathfrak{D}_H:= \{1 < |\zeta| < \rho_H\},
$$
such that the point $A_H$ is mapped to some point on the circle $\{|\zeta|=1\}$, say $\zeta=1$,  and $\varphi_H(B_H)$ is on the circle  $\{|\zeta|=\rho\}$.
We have
$$
m(\Omega_H)=m(\mathfrak{D}_H) = \frac{1}{2\pi}\log\rho_H,
$$
therefore, by \eqref{momega}, $\lim\limits_{H\to\infty}\rho_H=1$.
\begin{lemma1}\label{modulOmega_PQ}
We have
\begin{align}\label{ctOmega_PQ}
m^{-1}(\Omega_H)\sim m(\mbox{\boldmath $P$}_H)+m(\mbox{\boldmath $Q$}_H),\quad\mbox{\rm as }\quad H\to+\infty.
\end{align}
\end{lemma1}
\begin{proof} To show this, we will establish that for every $\varepsilon > 0$ there exists $H_0$ such that for all $H \geq H_0$ the images of the segments $A_HB_H$ and $C_HD_H$ under the map $\varphi_H$ lie in the intersections of the annulus $\mathfrak{D}_H$ and circular sectors of solution less than $\varepsilon$ centered at the origin.

First we consider the images of the segment $A_HB_H$.
The modulus is a conformal invariant, consequently, it is invariant under the shifts of the domain. Therefore, instead of $\Omega$ we can consider the domain which is obtained from $\Omega$ by translating along the real axis by the vector $(-c, 0)$. For convenience, we also denote the obtained domain by $\Omega$. In this case, we really have $c=0$,  therefore  the points $A_H=A$, $B_H=B$ and  the segment $A_HB_H=AB$ do not depend on~$H$.

Denote $r_H=\frac{2\rho_H}{1+\rho_H^2}$. It is clear that  $r_H<1$ and  $r_H\to1$ as $H\to \infty$.

The M\"obius transformation
$$
\psi_H(\zeta)=\rho_H\,\dfrac{\rho_H\zeta+i}{\zeta+i\rho_H}
$$
maps conformally the non-concentric annulus
$$
\mathfrak{G}_H:=\{\zeta\in\mathbb{C}\,:\,|\zeta+ir_H/2|>r_H/2,\, |\zeta|< 1\}.
$$
onto the annulus $\mathfrak{D}_H$.

Consider two cases.

a) Let $a<c=0$. Denote
$$
\mathfrak{G}:=\{\zeta\in\mathbb{C}\,:\,|\zeta+i/2|>1/2\},\, |\zeta|< 1\},
$$
$$
\Omega=\{(x,y):y<f_1(0)\}\setminus \{(x,f_2(0)):x\ge 0\}.
$$
By the definition of the kernel convergence, the family of domains $\Omega_H$ converges to $\Omega$ as to a kernel and also $\mathfrak{G}_H$ converges to $\mathfrak{G}$, as $H\to+\infty$.

 Now consider the conformal mappings  $\mathcal{F}_H:=\psi_H^{-1}\circ \varphi_H$. Then, by the Carath\'eodory convergence theorem~\cite[Chap.~5, Sect.~5]{goluzin1969geometric}, $\mathcal{F}_H$ converge locally uniformly in $\Omega$ to some conformal mapping $\mathcal{F}:\Omega\to\mathfrak{G}$, as $H\to+\infty$.

 Now with the help of Theorem~\ref{thr_rado_ge1}, we will show that the mappings $\mathcal{F}_H$ converge uniformly to the conformal mapping $\mathcal{F}:\Omega\to\mathfrak{G}$ on the segment $AB$. We should note that the limiting domain  $\Omega$ is not Jordan since it is a half plane cut along a ray and formally we can not employ Theorem~\ref{thr_rado_ge1} directly to the sequence $\mathcal{F}_H$ converging to $\mathcal{F}$. But we can apply to $\Omega_H$ additional conformal mappings
 $$
 g_H(z)=\delta_H\sqrt{\frac{D_H(z-A)}{D_H-z}}+c_H,
 $$
with appropriate fixed branches of the square root, where $\delta_H\to 1$ and $c_H\to0$, as $H\to\infty$.
Then the domains $g_H(\Omega_H)$ converge to a Jordan domain which is the image of $\Omega$ under the mapping $g(z)=\sqrt{z-A}$. Obviously, $g_H$ converges to $g$ uniformly on $AB$. We can choose  $\delta_H$ and $c_H$  such that $g_H(AB)\subset \overline{g(\Omega)}$. Let us fix a compact subset $K$ of  $g(\Omega)$ such that $g_H(AB)\subset \overline{g(\Omega)}$. By Theorem~\ref{thr_rado_ge1}, taking into account Remark~\ref{remfaminf}, we obtain that $\mathcal{F}_H\circ g_H^{-1}$ converge to $\mathcal{F}\circ g^{-1}$ on $K$. From this fact we conclude that $\mathcal{F}_H$ converge to $\mathcal{F}$ on $AB$.

Denote $\gamma_H=\varphi_H(AB)$.
Let us fix a neighborhood $U$ of the curve $\Gamma=\mathcal{F}(AB)$ which lies outside of some neighborhood $V$ of the point $(-i)$. Because of the uniform convergence, for sufficiently large $H$ the curves $\Gamma_H=\mathcal{F}_H(AB)$ lie in $U$. Consequently, they do not intersect $V$. We also note that for every $\delta>0$ the M{\"o}bius transforms $\psi_H$
converge uniformly to the constant $1$ on the set $|\zeta+i|\geq\delta$, as $H\to\infty$. This follows that for sufficiently large $H$ the curves $\gamma_H=\psi_H(\Gamma_H)$ lie in an arbitrarily small neighborhood of the point $1$.\smallskip

b) Let now $a=c$. This case is considered similarly to the case $a>c$; we should note that now the kernel of the sequence $\Omega_H$ equals to $\widetilde{\Omega}$ which is the complex plane with excluded rays
$\{(x,f_k(0)):x\ge 0\}$, $k=1$, $2$.

This means that for every $\varepsilon > 0$ there exists $H_0$ and $\alpha$ such that for all $H \geq H_0$ the image of the segments $A_HB_H$ under the map $\varphi_H$ lie in the sector $\alpha<\arg z<\alpha+\varepsilon$. Similarly we prove that for every $\varepsilon > 0$ there exists $H_1$ and $\beta $ such that for all $H \geq H_1$ the image of the segments $C_HD_H$ under the map $\varphi_H$ lie in the sector $\beta<\arg z<\beta +\varepsilon$. Without loss of generality we can assume that $H_0=H_1$ and $\alpha=0<\varepsilon<\beta<\beta+\varepsilon<2\pi$. Comparing the modules of $\varphi_H(\mbox{\boldmath$P$}_H)$ and $\varphi_H(\mbox{\boldmath$Q$}_H)$ with the modules of quadrilaterals lying in the annulus $\mathfrak{D}_H$ inside of the sectors bounded with the rays $\arg z = 0, \varepsilon, \beta,$ and $\beta+\varepsilon$, taking into account the conformal invariance of moduli, we conclude that for $H\ge H_0$
$$
(\beta-\varepsilon)/\log \rho_H<m(\mbox{\boldmath$Q$}_H)<(\beta+\varepsilon)/\log \rho_H,
$$
$$
(2\pi-\beta-\varepsilon)/\log \rho_H<m(\mbox{\boldmath$P$}_H)<(2\pi-\beta+\varepsilon)/\log \rho_H,
$$
Adding the above inequalities, we obtain
\begin{align}\label{modul_PH_QH_OmegaH}
(2\pi-2\varepsilon)/\log \rho_H<m(\mbox{\boldmath$Q$}_H)+m(\mbox{\boldmath$P$}_H)<(2\pi+2\varepsilon)/\log \rho_H, \quad\text{ for all }\quad H\ge H_0.
\end{align}
Because, $\varepsilon>0$ can be chosen arbitrarily small, from \eqref{modul_PH_QH_OmegaH} we get
$$
m(\mbox{\boldmath$Q$}_H)+m(\mbox{\boldmath$P$}_H)\sim 2\pi/\log\rho_H=m^{-1}(\Omega_H).
$$

The lemma is entirely proven.
\end{proof}

Now we  will estimate the behavior of the conformal modulus of $\mbox{\boldmath$P$}_H$ as $H\to+\infty$.

\begin{lemma1}\label{mdmodulP_H}
We have
\begin{align}\label{ctmodulP_H}
m(\mbox{\boldmath$P$}_H)=O(\log H),\quad \mbox{\rm as }\quad H\to+\infty.
\end{align}
\end{lemma1}
\begin{proof}
As in the proof of Lemma~\ref{modulOmega_PQ}, we shift $\Omega$ by the vector $(-c, 0)$ and keep for it and its subdomains the same notation.

Consider the quadrilateral $\mbox{\boldmath$P$}$ (Fig.~\ref{f6}).  Let us fix some real numbers $\kappa, \ell$ and positive $\varepsilon>0$ such that $f_2(c)<\kappa-\varepsilon<\kappa< f_1(c)$, $f_2(d)<\ell-\varepsilon<\ell< f_1(d)$.
If we replace the vertices of $\mbox{\boldmath$P$}$ with the points $d+i\ell$, $d+i(\ell-\varepsilon)$, $c+i(\kappa-\varepsilon)$, $c+i\kappa$,
then the modulus of the obtained quadrilateral $\mbox{\boldmath$P$}^{1}$ will be greater than the modulus of $\mbox{\boldmath$P$}$. Now we choose some rays $L_1$ and $L_2$ which come from the points $c+i\kappa$ and $d+i\ell$ and lie in $\mathcal{P}$. If we replace the side of the quadrilateral $\mbox{\boldmath$P$}^{1}$ with end points   $c+i\kappa$ and $d+i\ell$ by the union of the rays $L_1$ and $L_2$, then we obtain a quadrilateral $\mbox{\boldmath$P$}^{2}$, the modulus of which is greater than $m(\mbox{\boldmath$P$}^{1})$. Further, for sufficiently large $M$ the segment with enpoints $c+i(\kappa-M)$ and  $d+i(\kappa-M)$ is lower than the graph of $y=f_2(x)$. Replacing the boundary arc of the quadrilateral  $\mbox{\boldmath$P$}^{2}$ connecting the points $c+i\kappa$ and $d+i\ell$ and containing the graph of $y=f_2(x)$ with the three-link polygonal line $L$ with vertices at the points $c+i\kappa$, $c+i(\kappa-M)$, $d+i(\ell-M)$  and $d+i\ell$, we obtain a quadrilateral $\mbox{\boldmath$P$}^{3}$, the modulus of which is greater than  $m(\mbox{\boldmath$P$}^{2})$. We point out that the boundary of the quadrilateral consists of the rays
$L_1$, $L_2$ and the polygonal line $L$; its vertices are $d+i\ell$, $d+i(\ell-\varepsilon)$, $c+i(\kappa-\varepsilon)$, and $c+i\kappa$.

Denote by $\mbox{\boldmath$P$}^{3}_H$ the quadrilateral which is obtained from $\mbox{\boldmath$P$}^{3}$ by the stretching with coefficient $H$. It is clear that
\begin{equation}\label{mphph3}
m(\mbox{\boldmath$P$}_H)\le m(\mbox{\boldmath$P$}^3_H),
\end{equation}
therefore, to prove the lemma, we only need to prove that $m(\mbox{\boldmath$P$}^{3}_H)= O(\log H)$, $H\to+\infty$.
Let the equations of the straight lines containing $L_1$, the segment with endpoints $c+i(\kappa-M)$, $d+i(\kappa-M)$ and $L_2$ have the form $y=a_kx+b_k $, $k=1$, $2$,~$3$. Then their images under the stretching are on the lines $y=(a_k/H)x+b_k$, $1\le k\le 3$. Now we define the mapping
$$
\eta_H:x+iy\mapsto x+iv(x,y),
$$
where
\begin{align}\label{defquasi1}
v(x,y)= {\begin{cases}
y-\frac{a_1}{H}\,x-b_1,\quad \mbox{\rm if}\quad x\leq Hc,\\[0.1cm]
y- \frac{a_2}{H}\,x-b_2 -M,\quad \mbox{\rm if }\quad Hc\leq x\leq Hd,\\[0.1cm]
y- \frac{a_3}{H}\,x-b_3,\quad \mbox{\rm if}\quad x\geq Hd.
\end{cases}}
\end{align}
It is a $K$-quasiconformal homeomorphism of the plane with
$$
K=K(H)=\frac{1+k}{1-k},\quad  k=\frac{\mathfrak{a}}{\sqrt{\mathfrak{a}^2+4H^2}},\quad \mathfrak{a}:=\max_{1\le k
\le 3} |a_k|.
$$
Since $K(H)\to 1$, $H\to+\infty$, there exists $K_0\geq 1$ such that $K(H)\le K_0$.

Denote by $\mbox{\boldmath$T$}_H$ the quadrilateral which is the image of $\mbox{\boldmath$P$}^{3}_H$ under the mapping $\eta_H$. Then, because of quasiinvariance of conformal modulus under quasiconformal mappings, we have
\begin{equation}\label{p3hk0}
m(\mbox{\boldmath$P$}^{3}_H)\le K(H) m(\mbox{\boldmath$T$}_H)\le  K_0 m(\mbox{\boldmath$T$}_H).
\end{equation}

Because of \eqref{mphph3} and  \eqref{p3hk0},  to prove Lemma~\ref{mdmodulP_H},
we only need to establish that
\begin{equation}\label{thoh}
m(\mbox{\boldmath$T$}_H)=O(\log H),\quad\mbox{\rm as}\quad H\to+\infty.
\end{equation}
Since conformal modulus is invariant under shifts, we can assume that $c=-d$. We note that $\mbox{\boldmath$T$}_H=(T_H; dH, dH-i\varepsilon, -dH-i\varepsilon, -dH)$ where $T_H$ is obtained from the lower half plane by removing the rectangle $[-dH, dH]\times[-M,0]$. The quadrilateral $\mbox{\boldmath$T$}_H$ is symmetric with respect to the imaginary axis. Its modulus is twice times greater than the modulus of  its right half, i.e. the quadrilateral $\mbox{\boldmath$T$}_H^1$ which is the fourth quarter of the plane with the excluded rectangle $[0, dH]\times[-M,0]$, with vertices  $dH$, $dH-i\varepsilon$, $-iM$, and $\infty$. If we replace the side of $\mbox{\boldmath$T$}_H^1$ which is the vertical ray, going down from the point $-iM$ to infinity, with the ray going from this point to the left and shift the obtained quadrilateral by the vector $-dH$ then we obtain the quadrilateral $\mbox{\boldmath$U$}_H=(U; 0, -i\varepsilon,-d_H-iM, \infty)$ the modulus of which is bigger than the modulus of $\mbox{\boldmath$T$}_H^1$. Here $U$ is the domain the boundary of which is the polygonal line consisting of the horizontal rays $\{(x, -M), x\le 0\}$, $\{(x, 0), x\ge 0\}$ and the vertical segment $[-iM, 0]$. We have
\begin{equation}\label{mthuh}
m(\mbox{\boldmath$T$}_H)\le 2m(\mbox{\boldmath$U$}_H).
\end{equation}

Now we need the following auxiliary result.

\begin{lemma1}\label{mdmodulU_H} We have
$$
m(\mbox{\boldmath$U$}_H)=O(\log H), \quad\mbox{\rm as }\quad H\to+\infty.
$$
\end{lemma1}

We note that this lemma, in fact,  is a generalization of  \cite[Sect.~4, Lemma~4.6]{nvgiang2021}.

\begin{proof} We can map conformally the lower half plane $\mathbb{L}$ onto $U$ by a function $g$ such that $g(-1) = -iM$, $g(1) = 0$, $g(\infty) = \infty$. The function $g$ has the form
\begin{equation}\label{gz}
g(\zeta)=\frac{M}{\pi}\int_1^\zeta\sqrt{\frac{\omega+1}{\omega-1}}\,d\omega
\end{equation}
where the branch of the square root is fixed such that it is positive for real $\omega>1$.
Denote by $(-R_H)$ and  $\xi$ the preimages of the points $(-dH-iM)$ and $(-i\varepsilon)$ under the mapping $g(\zeta)$, $R_H>0$, $-1<\xi<1$. From \eqref{gz} it follows that
$$
g(\zeta)\sim \frac{M}{\pi}\,\zeta,\quad\mbox{\rm as }\quad \zeta\to\infty,
$$
therefore,
$$
R_H\sim \frac{M}{\pi}\,dH.
$$
Further, the conformal modulus of $\mbox{\boldmath$U$}_H$ is equal to the conformal modulus of the quadrilateral $\mbox{\boldmath$L$}_H=(\mathbb{L}; \infty, -R_H, \xi, 0)$. In its turn, by the Riemann-Schwarz symmetry principle, $m(\mbox{\boldmath$L$}_H)$ is a half of the modulus of the doubly-connected domain $\mathcal{A}=\overline{\mathbb{C}}\setminus\left( (-\infty, -R_H]\cup [\xi, 1]\right)$, which is an image of the Teichm\"uller ring under a linear transform. Thus, we have (see, e.g., \cite[Chap.3]{ahlfors2006lectures}; \cite{avv})
$$
2\pi m(\mathcal{A})\sim \log R_H\sim\log H, \quad\mbox{\rm as }\quad H\to \infty.
$$
Since $m(\mbox{\boldmath$U$}_H)\le 2 m(\mathcal{A})$, we proved the statement of the lemma.
\end{proof}

From Lemma~\ref{mdmodulU_H} and \eqref{mthuh} it follows that \eqref{thoh} holds. This proves Lemma~\ref{mdmodulP_H}.

At last, combining \eqref{minv}, \eqref{ctOmega_PQ} and \eqref{ctmodulP_H} we obtain \eqref{ctulmd1}.
Therefore, Theorem~\ref{main} is completely proved.
\end{proof}

\section{Discussion}
There are still many interesting issues concerning investigation
of the behavior of conformal moduli of domains
under their geometric transformations. Our next task includes studying the behavior of the exterior conformal moduli of arbitrary quadrilaterals when stretched in the direction of the abscissa axis.
\section*{FUNDING}
The work of the second author is performed under the development program of the Volga Region Mathematical Center (agreement no.~075-02-2022-882).


\begin{thebibliography}{99}
\bibitem{ahlfors2006lectures}
Ahlfors~L.~V., \emph{Lectures on Quasiconformal Mappings, Vol.~38} (American Mathematical Society, 2006).

\bibitem{dautova2018}
D.~N.~Dautova and S.~R.~Nasyrov, ``Asymptotics of the modules of
mirror symmetric doubly connected domains under stretching,'' Mathematical Notes \textbf{103} (3), 537--549 (2018).

\bibitem{dautova2019}
D.~N.~Dautova and S.~R.~Nasyrov, ``Asymptotics of conformal module
of nonsymmetric doubly connected domain under unbounded stretching
along the real axis,'' Lobachevskii Journal
  of Mathematics \textbf{40} (9), 1268--1274 (2019).

\bibitem{dautova2021conformal}
Dautova~D., Nasyrov~S. and Vuorinen~M., ``Conformal Modulus of the Exterior of Two Rectilinear Slits,'' Computational Methods and Function Theory (Springer) \textbf{21} (1), 109--130 (2021).

\bibitem{avv}
Glen~D.~Anderson, Mavina~K.~Vamanamurthy, Matti~K.~Vuorinen, \emph{Conformal invariants, inequalities and quasiconformal maps} (Canadian Mathematical Society Series of Monographs and Advanced Texts. A Wiley-Interscience Publication. J. Wiley, 1997).

\bibitem{goluzin1969geometric}
G.~M.~Goluzin, \emph{Geometric Theory of Functions of a Complex
Variable} (American Mathematical Soc.,  Providence, RI,
Translations of Mathematical Monographs, Vol.~26, 1969).

\bibitem{hakula2021conformal}
 Hakula~H., Nasyrov~S. and Vuorinen~M., ``Conformal moduli of symmetric circular quadrilaterals with cusps,'' Electronic Transactions on Numerical Analysis (Kent State University) \textbf{54}, 460--482 (2021).

\bibitem{hrv1}
Hakula~H., Rasila~A., and Vuorinen~M., ``On moduli of rings and quadrilaterals: algorithms and experiments,'' SIAM Journal on Scientific Computing \textbf{33} (1), 279--302 (2011).

\bibitem{hrv2}
Hakula~H., Rasila~A., and Vuorinen~M., ``Computation of exterior moduli of quadrilaterals,'' Electron. Trans. Numer. Anal. \textbf{40}, 436--451 (2013).

\bibitem{lang2003complex}
Lang S., \emph{Complex analysis, Vol.~103} (Springer Science \& Business Media, 2003).

\bibitem{nass}
Nasser~Mohamed~MS, ``PlgCirMap: A MATLAB toolbox for computing conformal mappings from polygonal multiply connected domains onto circular domains,'' SoftwareX \textbf{11}, 100464 (2020).

\bibitem{nasser2021circular}
Nasser~M., Rainio~O., Rasila~A., Vuorinen~M., Wallace~T., Yu~H. and Zhang~X., ``Circular arc polygons, numerical conformal mappings, and moduli of quadrilaterals,'' arXiv preprint arXiv:2107.11485 (2021).

\bibitem{nasyrov2015riemann}
Nasyrov~S.~R., ``Riemann--{S}chwarz reflection principle and
asymptotics of modules of rectangular frames,'' Computational
Methods and Function Theory \textbf{15} (1), 59--74 (2015).

\bibitem{nvgiang2021}
Nasyrov~S.~R., Van Giang~N., ``Asymptotics of the conformal modulus of unbounded symmetric doubly-connected domain under stretching,'' Lobachevskii Journal
of Mathematics \textbf{42} (12), 2895--2904 (2021).

\bibitem{nasyrov2021moduli}
Nasyrov S.~R., Sugawa~T. and Vuorinen~M., ``Moduli of quadrilaterals and quasiconformal reflection,'' arXiv:2111.08304 (2021).

\bibitem{papamichael2010numerical}
Papamichael~N., Stylianopoulos~N., \emph{Numerical Conformal
Mapping: Domain Decomposition and the Mapping of Quadrilaterals}
(World Scientific, Singapore, 2010).

\bibitem{kuhnau2005conformal}
R.~K{\"u}hnau, ``The conformal module of quadrilaterals and of
rings,'' in: \emph{Handbook of Complex Analysis}, Ed. by
R.~K{\"u}hnau (Elsevier, 2005), Vol.~2, pp.~99--129.

\bibitem{remmert2013classical}
Remmert~R., \emph{Classical topics in complex function theory, Vol.~172} (Springer Science \& Business Media, 2013).

\bibitem{schiff1993normal}
Schiff~J.~L., \emph{Normal families} (Springer Science \& Business Media, 1993).

\bibitem{vuorinen2013exterior}
Vuorinen~M., Zhang~X., ``On exterior moduli of quadrilaterals
and special functions,'' Journal of Fixed Point Theory and
Applications \textbf{13} (1), 215--230 (2013).
\end{thebibliography}
\end{document}